# Analysis.

## An homotopy of isometries related to a probability density.

### Roland Groux

*Lycée Polyvalent Rouvière, rue Sainte Claire Deville,
BP 1205. 83070 Toulon .Cedex. France.*

Email : roland.groux@orange.fr

**Abstract.**

We are studying here a family of probability density functions indexed by a real parameter, and constructed from homographic relations between associated Stieltjes transforms. From the analysis of orthogonal polynomials we deduce a family of isometries in relation to the classical operators creating secondary polynomials and we give an application to the explicit resolution of specific integral equations.

## 1. Introduction and notations.

Let $\rho$ a probability density function on an interval $I$ bounded with $a$ and $b$. We note $c_n = \int_a^b x^n \rho(x) dx$ the moment of order $n$.

The Stieltjes transform of the measure of density $\rho$ is defined on $\mathbb{C} - I$ by the formula: $z \mapsto S_\rho(z) = \int_a^b \frac{\rho(t) dt}{z - t}$. (See [1], [2]).

We note $n \mapsto P_n$ a Hilbert base of normalized polynomials for the classic inner product: $(f, g) \mapsto < f / g >_\rho = \int_a^b f(t) g(t) \rho(t) dt$ on the associated Hilbert' space $L^2(I, \rho)$.

$T_\rho$ denotes the operator $f(x) \mapsto g(x) = \int_a^b \frac{f(t) - f(x)}{t - x} \rho(t) dt$ creating secondary polynomials $Q_n = T_\rho(P_n)$.

Let us recall the results below: (See [3])

Let be a positive measure on $I$ associated to a density function $\mu$, also allowing moments of any order and having Stietjes's transformation linked to these of $\rho$ by the equality: $S_\mu(z) = z - c_1 - \frac{1}{S_\rho(z)}$. We can then conclude:

_ Secondary polynomials $A_n = Q_{n+1}$ relative with $\rho$ then form an orthonormal family for the inner product induced by $\mu$.

_ Secondary polynomials $B_n = T_\mu(A_n)$ relative with $\mu$ are defined by the formula: $B_n(x) = (x - c_1) Q_{n+1}(x) - P_{n+1}(x)$

_ If the density $\rho$ is continuous and provided the existence of $\varphi(x) = \lim_{\varepsilon \to 0^+} 2\int_a^b \frac{(x-t)\rho(t)dt}{(x-t)^2 + \varepsilon^2}$, we can make $\mu$ explicit by the formula: $\mu(x) = \frac{\rho(x)}{\frac{\varphi^2(x)}{4} + \pi^2 \rho^2(x)}$.

_ The operator $f(x) \mapsto g(x) = \int_a^b \frac{f(t) - f(x)}{t - x} \rho(t)dt$ creating secondary polynomials extends to a continuous linear map linking the space $L^2(I, \rho)$ to the Hilbert'space $L^2(I, \mu)$, whose restriction to the hyperplane $H_\rho$ of the functions orthogonal for $\rho$ with $P_0 = 1$ constitutes an isometric function for both norms respectively.

Under the mentioned assumptions, the function $\varphi$ presented above will be call the *reducer* of $\rho$, the measure of density $\mu$ will be call *secondary measure* associated with $\rho$. Its moment of order 0 is equal to $d_0 = c_2 - (c_1)^2$. If we normalize $\mu$, we introduce $\mu_0 = \frac{\mu}{d_0}$ a probability density function called the '*normalized secondary measure*' of $\rho$.

By definition we said in what follows that two measures are equi-normal if they lead to the same normalized secondary measure.

## 2. A family of equi-normal measures.

In this section we suppose that $I$ is a compact interval.

We consider in this section a real parameter $t > 0$ and a density of probability $\rho_t$ whose secondary measure density is $t\mu$. We also requires that the moment of order 1 of the two densities $\rho$ and $\rho_t$ are the same. Under these assumptions, the coupling of Stieltjes tansforms results in:

$S_\mu(z) = z - c_1 - \frac{1}{S_\rho(z)}$ and $S_{t\mu}(z) = tS_\mu(z) = z - c_1 - \frac{1}{S_{\rho_t}(z)}$.

We deduce the Stieltjes tansform of $\rho_t$: $S_{\rho_t}(z) = \frac{S_\rho(z)}{t + (1-t)(z - c_1)S_\rho(z)}$. (2.1)

We recall the formula of Stieltjes-Perron making explicit the density from its Stieltjes transform: $\rho(x) = \lim_{\varepsilon \to 0^+} \frac{S_\rho(x - i\varepsilon) - S_\rho(x + i\varepsilon)}{2i\pi}$ and the definition of the reducer of $\rho$: $\varphi(x) = \lim_{\varepsilon \to 0^+} S_\rho(x + i\varepsilon) + S_\rho(x + i\varepsilon)$.

We easily deduce: $\lim_{\varepsilon \to 0^+} S_\rho(x + i\varepsilon).S_\rho(x + i\varepsilon) = \frac{\varphi^2(x)}{4} + \pi^2 \rho^2(x)$.

By applying Stieltjes-Perron to $S_{\rho_t}$, we obtain after simplifying:

$$\rho_t(x) = \frac{t\rho(x)}{[(t-1)(x-c_1)\frac{\varphi(x)}{2} - t]^2 + \pi^2\rho^2(x)(t-1)^2(x-c_1)^2} \quad (2.2)$$

However, as shown by the following examples, the function defined by the above formula does not always its Stieltjes transform equals to $H(z) = \frac{S_\rho(z)}{t + (1-t)(z-c_1)S_\rho(z)}$.

It will be easily seen by noting that the denominator of $H(z)$ may vanish outside $I$ if the parameter $t$ becomes too high.

But if effectively this function $H$ is the Stieltjes transform of $\rho_t$, then we deduce from $S_\mu(z) = z - c_1 - \frac{1}{S_\rho(z)}$ the formula: $S_{\rho_t}(z) = \frac{1}{z - c_1 - tS_\mu(z)}$.

So the moment of order 0 of $\rho_t$ is obtained by : $\lim_{|z|\to\infty} zS_{\rho_t}(z) = \lim_{|z|\to\infty} \frac{1}{1 - \frac{c_1}{z} - \frac{tS_\mu(z)}{z}} = 1$.

Thus, $\rho_t$ is a probability density function and its moment of order 1 is $c_1$ because:

$$\lim_{|z|\to\infty} z^2(S_{\rho_t}(z) - \frac{1}{z}) = \lim_{|z|\to\infty} \frac{c_1 - tS_\mu(z)}{1 - \frac{c_1}{z} - \frac{tS_\mu(z)}{z}} = c_1.$$

We can finally conclude thanks to $tS_\mu(z) = z - c_1 - \frac{1}{S_{\rho_t}(z)}$, that the secondary measure of the measure of density $\rho_t$ is $t\mu$. So $\rho_t$ is effectively equi-normal with $\rho$.

Now consider some examples. We will study particularly the function
$$t \mapsto f(t) = t\int_I \frac{\rho(x)dx}{[(t-1)(x-c_1)\frac{\varphi(x)}{2} - t]^2 + \pi^2\rho^2(x)(t-1)^2(x-c_1)^2}$$ making explicit the moment
of order 0 for $\rho_t$.

- The Tchebychev measure of the second kind over $]-1,1[$

$\rho(x) = \frac{2}{\pi}\sqrt{1-x^2}$ , and we have $c_1 = 0$ ; $\varphi(x) = 4x$ and $\mu(x) = \frac{\rho(x)}{4}$

The function $f$ is made explicit here by: $t \mapsto f(t) = \frac{2t}{\pi}\int_{-1}^{1} \frac{\sqrt{1-x^2}}{t^2 + 4(1-t)x^2}dx$

A quick viewing using MAPLE shows that *f* keeps the constant value 1 all over the interval ]0,2].

```
> f:=proc(t)
local x,g;
g:=evalf((2*t/Pi)*int(sqrt(1-x^2)/(t^2+4*(1-t)*x^2),x=-1..1));end :
```

```
> plot(f);
```

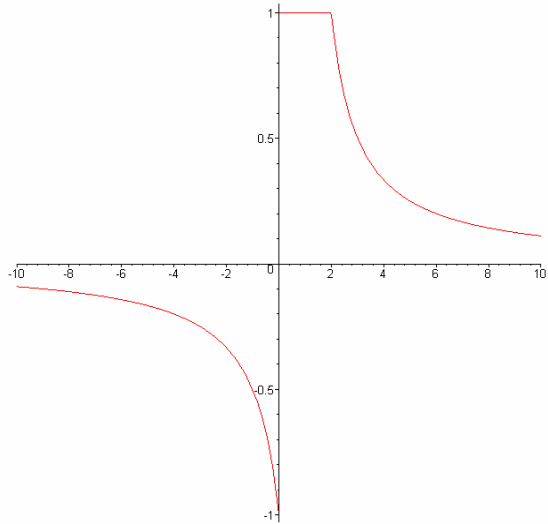

This can be directly verified by direct calculation, as well as the following formulas:

$$\rho_t(x) = \frac{2t\sqrt{1-x^2}}{\pi[t^2 + 4(1-t)x^2]}, \text{ with } t \text{ in } ]0, 2].$$

The secondary measure of the density $\rho_t$ is $\mu_t = \frac{t}{4}(\frac{2\sqrt{1-x^2}}{\pi})$

The reducer of $\rho_t$ is defined by : $\varphi_t(x) = \frac{2(4-2t)x}{t^2 + 4(1-t)x^2}$.

Stieltjes transform is given by: $S_t(z) = \frac{2}{(2-t)z + t\sqrt{z^2 - 1}}$

For *t*=1 we find of course the Tchebychev measure of the second kind.

For *t*=2 we obtain the Tchebychev measure of the first kind.

For $t = \frac{4}{3}$ we have $\rho(x) = \frac{6\sqrt{1-x^2}}{\pi(4-3x^2)}$ and $\varphi(x) = \frac{6x}{4-3x^2}$.

- **The uniform Lebesgue's measure over [0,1].**

Uniform density $\rho(x) = 1$ with the reducer $\varphi(x) = 2\ln(\frac{x}{1-x})$ and $c_1 = \frac{1}{2}$.

Here we have : $t \mapsto f(t) = t\int_0^1 \frac{dx}{[(t-1)(x-\frac{1}{2})\ln(\frac{x}{1-x}) - t]^2 + \pi^2(t-1)^2(x-\frac{1}{2})^2}$

> `f:=proc(t)`
> `local x,g;`
> `g:=evalf(t*int(1/(((t-1)*(x-1/2)*ln(x/(1-x))-t)^2+Pi*Pi*(t-1)^2*(x-1/2)^2),x=0..1));end :`
> `f(1.3);`
>                                   0.9799849175
>
> `plot(f);`

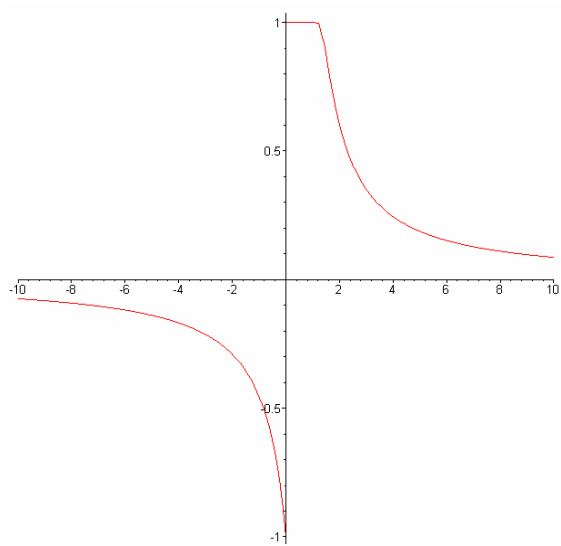

So it seems that, for any $t$ in $]0,1]$, the secondary measure of the density

$\rho_t(x) = \frac{t}{[(t-1)(x-\frac{1}{2})\ln(\frac{x}{1-x}) - t]^2 + \pi^2(t-1)^2(x-\frac{1}{2})^2}$ is : $\mu_t(x) = \frac{t}{[\ln^2(\frac{x}{1-x}) + \pi^2]}$.

We can obtain by MAPLE verification of this proportionality by examining the moments of order two.

Recall the formulas : $S_\mu(z) = z - c_1 - \frac{1}{S_\rho(z)}$ et $tS_\mu(z) = z - c_1 - \frac{1}{S_{\rho_t}(z)}$ .

Thus : $c_2' - (c_1')^2 = t(c_2 - (c_1)^2)$.

Or in this case: $c_1 = c_1' = \frac{1}{2}$ and $c_2 = \frac{1}{3}$. So we get: $c_2' = \frac{t+3}{12}$, confirmed by MAPLE :

```
> f2:=proc(t)
local x,g;
g:=evalf(t*int(x^2/(((t-1)*(x-1/2)*ln(x/(1-x))-t)^2+Pi*Pi*(t-
1)^2*(x-1/2)^2),x=0..1));end :

> plot(f2,0..1);
```

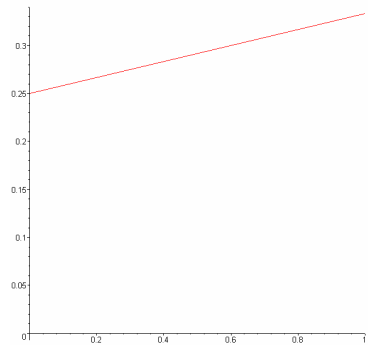

- For $\rho(x) = 2x$ over [0, 1], the reducer is $\varphi(x) = -4x\ln\left(\dfrac{1-x}{x}\right) - 4$ and $c_1 = \dfrac{2}{3}$.

```
> ro:=x->2*x:
> phi:=x->-4*(x*ln((1-x)/x)+1):
> c1:=evalf(int(x*ro(x),x=0..1)):
> f:=proc(t)
local x,g;
g:=evalf(t*int(ro(x)/(((t-1)*(x-c1)*phi(x)/2-t)^2+Pi^2*(t-
1)^2*ro(x)^2*(x-c1)^2),x=0..1));end :

> f(0.45);plot(f);
```
$$0.9999999999$$

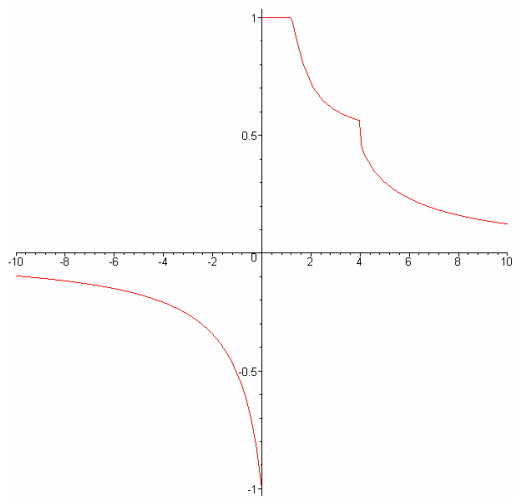

- For $\rho(x) = \dfrac{3\sqrt{x}}{2}$ over [0, 1], we have $\varphi(x) = 3\text{LerchPhi}(x,1,-\dfrac{1}{2})$ and $c_1 = \dfrac{3}{5}$.

```
> ro:=x->3*sqrt(x)/2 ;     ro := x → (3/2)√x

> c1:=int(x*ro(x),x=0..1) ;    c1 := 3/5

> phi:=x->3*LerchPhi(x,1,-1/2);   φ := x → 3 LerchPhi(x, 1, -1/2)

> f:=proc(t)
local x,g;
g:=evalf(int(ro(x)/(((t-1)*(x-c1)*phi(x)/2-
t)^2+Pi^2*ro(x)^2*(t-1)^2*(x-c1)^2),x=0..1));evalf(t*g);end ;
    f := proc(t)
      local x, g;
        g := evalf(int(ro(x)/(
            (1/2×(t − 1)×(x − c1)×φ(x) − t)^2 + π^2×ro(x)^2×(t − 1)^2×(x − c1)^2),
            x = 0 .. 1));
          evalf(t×g)
      end proc

> f(0.6);  1.000000000
> f(0.45628);  0.9999999999
> f(0.2157);  1.000000000
> f(2);  0.7496041742
> f(1.24);  0.9911159300
```

Now we'll try to explain why we get a good probability density equi-normal with $\rho$ when the parameter $t$ remains in the interval $]0, 1]$.

Consider first the following lemma: (For an easy writing we note in what follows $I = [-1, 1]$)

### 2.3. Lemma.

The equation $\boxed{t + (1-t)(z - c_1)S_\rho(z) = 0}$ has no solutions belonging to the open set $O = \mathbb{C} - I$ when the parameter $t$ is element of $]0, 1]$.

**Proof**.

Noting $m = \dfrac{t}{t-1}$ and $z = x + iy$, the equation becomes: $\displaystyle\int_I \dfrac{(x - c_1 + iy)\rho(u)du}{x - u + iy} = m$

Separating the real and imaginary we get:

$$\int_I \dfrac{x^2 + y^2 - u(x - c_1) - c_1 x}{(x-u)^2 + y^2}\rho(u)du = m \quad \text{and} \quad y\int_I \dfrac{c_1 - u}{(x-u)^2 + y^2}\rho(u)du = 0$$

Let us consider first the solutions out of the real axis.

If $y \neq 0$ we have from imaginary part above: $\int_I \dfrac{u}{(x-u)^2 + y^2} \rho(u)du = \int_I \dfrac{c_1}{(x-u)^2 + y^2} \rho(u)du$.

And for real part : $\int_I \dfrac{x^2 + y^2 + c_1^2 - 2c_1 x}{(x-u)^2 + y^2} \rho(u)du = m \int_I \dfrac{x^2 + y^2 + u^2 - 2c_1 x}{(x-u)^2 + y^2} \rho(u)du$

After reductions : $\int_I \dfrac{(1-m)[(x-c_1)^2 + y^2] + m(c_1^2 - u^2)}{(x-u)^2 + y^2} \rho(u)du = 0$

And with the variable $t$ : $\int_I \dfrac{([(x-c_1)^2 + y^2] + t(u^2 - c_1^2)}{(x-u)^2 + y^2} \rho(u)du = 0$

Thus, if $c_1 = 0$, there are no non-real solutions when $t > 0$.

Study now real solutions.

For $y = 0$, the equation easily simplifies to: $\int_I \dfrac{m.u + (1-m)x - c_1}{x - u} \rho(u)du = 0$.

And with the variable $t$ : $\int_I \dfrac{t.u - x + c_1(1-t)}{x - u} \rho(u)du = 0$

Here again if $c_1 = 0$, there are no real solutions out of $I$ if $0 < t < 1$.
However there may exist if $t > 1$ as shown in the following example for a Tchebychev's measure and the value $t = 3$.

```
> ro:=x->sqrt(1-x^2):
> f:=proc(x)
local t,f;
f:=evalf(int((3*t-x)*ro(t)/(x-t),t=-1..1));end;
```
$f := \mathbf{proc}(x) \mathbf{\ local\ } t, f; f := \text{evalf}(\text{int}((3 \times t - x) \times \text{ro}(t)/(x-t), t = -1 .. 1)) \mathbf{\ end\ proc}$

```
> plot(f);
```

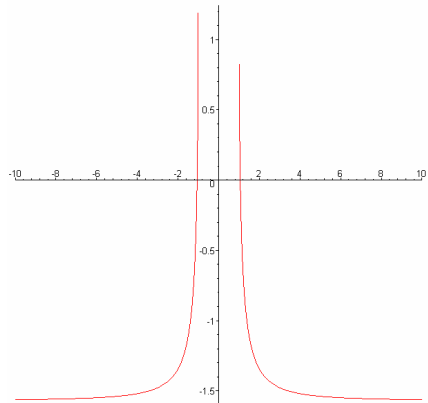

```
> f(1.06);  0.005893680455
> f(1.07); -0.07793247767
```

Now let us analyze the general case where : $c_1 \neq 0$ and $t \in ]0,1]$

> Resuming the first study of non-real solutions. We got the qualities:

$$\int_I \frac{u}{(x-u)^2 + y^2}\rho(u)du = \int_I \frac{c_1}{(x-u)^2 + y^2}\rho(u)du \text{ and } \int_I \frac{([(x-c_1)^2 + y^2] + t(u^2 - c_1^2))}{(x-u)^2 + y^2}\rho(u)du = 0.$$

Or $[(x-c_1)^2 + y^2] + t(u^2 - c_1^2) = x^2 - 2c_1x + y^2 + tu^2 + (1-t)c_1^2$

The first integral equality leads to : $\int_I \frac{-2xu}{(x-u)^2 + y^2}\rho(u)du = \int_I \frac{-2xc_1}{(x-u)^2 + y^2}\rho(u)du$

(This by multiply by -2x). So we deduce the transform of the second integral:

$$J = \int_I \frac{([(x-c_1)^2 + y^2] + t(u^2 - c_1^2))}{(x-u)^2 + y^2}\rho(u)du = 0 = \int_I \frac{([(x-u)^2 + y^2] + (t-1)(u^2 - c_1^2))}{(x-u)^2 + y^2}\rho(u)du = K$$

Now write : $0 = (1-t)J + tK$. After simplifications this results to:

$$\int_I \frac{((1-t)(x-c_1)^2 + y^2 + t(x-u)^2)}{(x-u)^2 + y^2}\rho(u)du = 0 \text{ , impossible, because } t \in ]0,1]$$

> Now examine real solutions. We saw earlier in this case :

$\int_I \frac{t.u - x + c_1(1-t)}{x-u}\rho(u)du = 0$. Or we can write denominator as : $-t(x-u) + (1-t)(c_1 - x)$.

From this we deduce: $\int_I \frac{(1-t)(c_1 - x)}{x-u}\rho(u)du = t$ .

If we write $c_1 - x = c_1 - u + u - x$ , we gets : $t - 1 + \int_I \frac{(1-t)(c_1 - u)}{x-u}\rho(u)du = t$

So : $\int_I \frac{(1-t)(c_1 - u)}{x-u}\rho(u)du = 1$

Or if $x$ is located 'right on the interval' $I$, therefore higher than any element $u$ of $I$ and also higher than $c_1$, we have : $\frac{c_1 - u}{x-u} < \frac{x-u}{x-u} = 1$.

So we deduce if $t < 1$ : $1 \leq \int_I \frac{(1-t)(x-u)}{x-u}\rho(u)du = 1 - t$ . So we conclude $t \leq 0$, that is contrary to the hypothesis.

(We obtain a similar contradiction if $x$ is located 'left on the interval' $I$.
This completes the proof of Lemma (2.3)

Now study the following result:

### 2.4. Lemma.

Let $H$ a function of complex variable, holomorphic on $O = \mathbb{C} - I$.

We note $\tilde{I} = [-a, a]$, with $a > 1$ and we suppose here that:

2.4.1 $|H(z)|$ has limit 0 when the module of $z$ tends to infinity.

2.4.2 The two sequences of function defined on $\tilde{I}$ by $x \mapsto H(x - \frac{i}{n})$ and $x \mapsto H(x + \frac{i}{n})$ have a limit for the classic norm of $L^1([-a, a])$ when the index $n$ tends to infinity.

Under these hypotheses we can conclude that $H$ is the Sieltjes transform of $v$ defined on $I$ by

$$v(x) = \lim_{\varepsilon \to 0^+} \frac{H(x - i\varepsilon) - H(x + i\varepsilon)}{2i\pi}.$$

**Proof**.

Let $Z$ a complex number element of $O = \mathbb{C} - I$.

The function : $z \mapsto f(z) = \frac{H(z)}{(z - Z)}$ is holomorphic on $O - \{Z\}$ an its residue at the point $Z$ is $H(Z)$. Let $\lambda$ a fixed real with $0 < \lambda < a$

We consider $R_\varepsilon$ a rectangle with sides parallel to the axes, oriented in the direct sense, and whose summits have for affix : $1 + \lambda + \varepsilon(\pm i); -1 - \lambda + \varepsilon(\pm i)$, with $\varepsilon > 0$.

We note $C_b$ the circle focused on the origin, oriented in the direct sense and whose radius equal $b$, with $b > |Z|$.

The two paths $R_\varepsilon$ et $C_b$ are homotopic in $O$. Thus, according with Cauchy's theorem:
$$\int_{R_\varepsilon} f(z) dz = \int_{C_b} f(z) dz - 2i\pi H(Z). \quad (2.4.3)$$

Now study the evolution of this equality when: $b \to +\infty$ and $\varepsilon \to 0^+$.

➢ It is clear that for the two segments parallel to the imaginary axis, the contribution of the first integral vanished to 0 with $\varepsilon$, because $f$ has a finite limit at $1 + \lambda$ and $-1 - \lambda$.
$$\lim_{\varepsilon \to 0^+} \int_{-\varepsilon}^{\varepsilon} f(1 + \lambda + it) i \, dt = 0 = \lim_{\varepsilon \to 0^+} \int_{-\varepsilon}^{\varepsilon} f(-1 - \lambda + it) i \, dt.$$

➢ For the two parallels to the real axis, we obtain by adding the integrals over reversed paths:

$$\lim_{\varepsilon \to 0^+} \int_{-1-\lambda}^{1+\lambda} \left( \frac{H(t - i\varepsilon)}{t - i\varepsilon - Z} - \frac{H(t + i\varepsilon)}{t + i\varepsilon - Z} \right) dt = \lim_{\varepsilon \to 0^+} \left( \int_{-1-\lambda}^{1+\lambda} \frac{(t - Z)[H(t - i\varepsilon) - H(t + i\varepsilon)] dt}{(t - Z - i\varepsilon)(t - Z + i\varepsilon)} + J(\varepsilon) \right)$$

with $J(\varepsilon) = i\varepsilon \int_{-1-\lambda}^{1+\lambda} \frac{H(t-i\varepsilon) + H(t+i\varepsilon)}{(t-Z-i\varepsilon)(t-Z+i\varepsilon)} dt$. From (2.4.2) we easily deduce : $\lim_{\varepsilon \to 0^+} J(\varepsilon) = 0$, and also : $\lim_{\varepsilon \to 0^+} \int_{-1-\lambda}^{1+\lambda} \frac{(t-Z)[H(t-i\varepsilon) - H(t+i\varepsilon)]dt}{(t-Z-i\varepsilon)(t-Z+i\varepsilon)} = 2i\pi \int_{-1-\lambda}^{1+\lambda} \frac{v(t)dt}{t-Z}$

➢ Note that according Jordan's lemma, and thanks to (2.4.1), the integral on the circle $C_b$ vanished to 0 when radius $b$ tends to infinity. Thus, by limit passages, the formula (2.4.3) becomes: $2i\pi \int_{-1-\lambda}^{1+\lambda} \frac{v(t)dt}{t-Z} = -2i\pi H(z)$.

When $\lambda$ vanished to 0 we finally obtain the expected formula: $\int_{-1}^{1} \frac{v(t)dt}{Z-t} = S_v(Z) = H(Z)$

Thanks to these two precedent results we can now explain the importance of placing the parameter $t$ in the interval $]0,1]$. So in what follows we suppose $0 < t \leq 1$.

We'll apply the lemma (2.4) to the function $z \mapsto H(z) = \frac{S_\rho(z)}{t + (1-t)(z-c_1)S_\rho(z)}$.

First note that thanks to lemma (2.3), the function $H$ is holomorphic on $O = \mathbb{C} - I$

(2.4.1) is provided here directly by $\lim_{|z| \to \infty} S_\rho(z) = 0$ and $\lim_{|z| \to \infty} zS_\rho(z) = 1$

If we suppose that the density $\rho$ is such that the function $H$ above satisfies (2.4.2) we can effectively conclude that the Stieltjes transform of $\rho_t$ is $H$. In this case we have a family of density of probability equi-normal with the initial density $\rho$ and indexed with the real parameter $t \in ]0,1]$.

We can notice that when $t$ tends to 0, the limit of this family within the sense of distributions is the Dirac's measure $\delta_{c_1}$ concentrated at the first moment $c_1$. Indeed we check easily that for any continuous function $g$ on the interval $I$: $\lim_{t \to 0} \int_I f(u)\rho_t(u)du = f(c_1)$

Note also that the reducer $\varphi_t(x)$ of $\rho_t(x)$ tends towards the Dirac's reducer : $\frac{2}{x-c_1}$.

(The Sieltjes transform of $\delta_{c_1}$ is $S_{\delta_{c_1}}(z) = \frac{1}{z-c_1}$.)

Always within the sense of distributions we note also : $\lim_{t \to 0^+}(x-c_1)^2 \frac{\rho_t}{t} = \mu$.

## 3. Orthogonal polynomials for the density $\rho_t$.

We note here $n \mapsto P_n^t$ a classic sequence of orthonormal polynomials related to the density $\rho_t$ and $Q_n^t = T_{\rho_t}(P_n^t)$ the associated secondary polynomials. In accordance with the notations of the first paragraph, we have: $P_n^1 = P_n$ and $Q_n^1 = Q_n$.

Recall the important result mentioned in this introduction section: (3.1)

_ Secondary polynomials $A_n = Q_{n+1}$ relative with $\rho$ then form an orthonormal family for the inner product induced by $\mu$.

_ Secondary polynomials $B_n = T_\mu(A_n)$ relative with $\mu$ are defined by the formula:
$B_n(x) = (x - c_1)Q_{n+1}(x) - P_{n+1}(x)$

When we normalize $\mu$ into $\mu_0 = \dfrac{\mu}{d_0}$, the set of primary orthonormal polynomials becomes $\tilde{A}_n = \sqrt{d_0} A_n$ and the associated secondary polynomials are now made explicit by $\tilde{B}_n = \dfrac{B_n}{\sqrt{d_0}}$ because $T_{\mu_0} = \dfrac{T_\mu}{d_0}$.

The same normalization from the density equi-normal $\rho_t$ leads to the same set of polynomials, because normalize secondary measure of the density $\rho_t$ is also $\mu_0$.

So, thanks to the formulas (3.1) we get:

$$\begin{cases} \sqrt{d_0} Q_{n+1}(x) = \sqrt{d_0 t}\, Q_{n+1}^t(x) \\ \dfrac{1}{\sqrt{d_0}}[(x-c_1)Q_{n+1}(x) - P_{n+1}(x)] = \dfrac{1}{\sqrt{d_0 t}}[(x-c_1)Q_{n+1}^t(x) - P_{n+1}^t(x)] \end{cases}$$

Then, for $n \geq 1$:

$$(3.2) \quad \begin{cases} P_n^t(x) = \dfrac{1}{\sqrt{t}}[tP_n(x) + (1-t)(x-c_1)Q_n(x)] \\ Q_n^t(x) = \dfrac{1}{\sqrt{t}} Q_n(x) \end{cases}$$

**Associated isometries.**

From $Q_n = T_\rho(P_n)$ and (3.2), we see that the map $f(x) \mapsto \dfrac{1}{\sqrt{t}}[tf(x) + (1-t)(x-c_1)T_\rho(f)(x)]$ transforms $(P_n)_{n \geq 1}$ orthonormal for the density $\rho$ into $(P_n^t)_{n \geq 1}$ orthonormal for the density $\rho_t$.

According to Cauchy's theorem, we can extends this application to an isometric map linking the hyperplane $H_\rho$ of the functions of $L^2(I,\rho)$ orthogonal for $\rho$ with $P_0 = 1$ to the hyperplane $H_\rho^t$ constituted by the elements of $L^2(I,\rho_t)$ orthogonal with $P_0^t = 1$.

If we introduce $\tilde{\rho}_t = \dfrac{\rho_t}{t}$, we deduce that $V_\rho^t$ défined by :

$f(x) \mapsto V_\rho^t(f(x)) = tf(x) + (1-t)(x-c_1)T_\rho(f(x))$ is an isometric function linking $H_\rho$ equipped with the norm of $L^2(I,\rho)$ to $H_\rho^t$ equipped with the norm associated to $\tilde{\rho}_t$. This translates into the following equal:

$$(3.3) \boxed{\int_I f^2(x)\rho(x)dx = \int_I [tf(x)+(1-t)(x-c_1)T_\rho(f)(x)]^2 \frac{\rho_t(x)dx}{t} = \int_I [V_\rho^t(f)(x)]^2 \tilde{\rho}_t(x)dx}$$

For example, for the Tchebychev's measure of the second kind over $]-1,1[$ : $\rho(x) = \dfrac{2}{\pi}\sqrt{1-x^2}$,

we have $\tilde{\rho}_t(x) = \dfrac{2\sqrt{1-x^2}}{\pi[t^2+4(1-t)x^2]}$ ; $c_1 = 0$. So (3.3) translates in this case to:

$$\boxed{\int_{-1}^1 f^2(x)\sqrt{1-x^2}\,dx = \int_{-1}^1 [tf(x)+(1-t)xT_\rho(f(x))]^2 \frac{\sqrt{1-x^2}\,dx}{t^2+4(1-t)x^2}}$$, under the assumption that

$\bar{f} = \dfrac{2}{\pi}\int_{-1}^1 f(x)\sqrt{1-x^2}\,dx$ equals 0. (The formula (3.3) requires belonging of $f$ to $H_\rho$). In the general case we adjust the formula by changing $f$ into $f - \bar{f}$, as shown in the following MAPLE program.

```
>rho:=x->2*sqrt(1-x^2)/Pi;
>rho2:=(x,s)->2*sqrt(1-x^2)/(Pi*(s^2+4*(1-s)*x^2)):
>V:=proc(f,t)
local u,x,g,y,C;
C:=evalf(int(f(u)*rho(u),u=-1..1));
g:=t*(f(x)-C)+(1-t)*x*int((f(u)-f(x))*rho(u)/(u-x),u=-1..1);
y:=evalf(int(g^2*rho2(x,t),x=-1..1));end :
>V(x->x^3-2/(x+5)+1/(x^2+3),1.35);
                0.1010020263902955321

>W:=proc(g)
local x,C,y;
C:=evalf(int(g(x)*rho(x),x=-1..1));
y:=evalf(int(g(x)*g(x)*rho(x),x=-1..1)-C^2);end:

>W(x->x^3-2/(x+5)+1/(x^2+3));
                0.1010020263902957
```

We will see now that under an hypothesis of density related to $\rho_t$, we have the composition pattern : $V_\rho^t = (T_\rho^t)^{-1} \circ T_\rho$ , with the following usual notations:

- $T_\rho$ is the isometric function linking $H_\rho$ with the norm of $L^2(I,\rho)$ to $L^2(I,\mu)$.

- $T_\rho^t$ linking $H_\rho^t$ with the norm of $L^2(I,\rho_t)$ to $L^2(I,t\mu)$.

The existence of $T_\rho^t$ requires the density of the space of polynomials in $L^2(I,\rho_t)$, which is verified for instance if $I$ is a compact interval. Under this hypothesis, and thanks to the density of polynomials in $H_\rho$, it is sufficient to study the transform of $P_n$.

By definitions: $V_\rho^t(P_n) = tP_n + (1-t)(x-c_1)Q_n$

Thanks (3.2) and linearity of $T_\rho^t$, we have: $T_\rho^t \circ V_\rho^t(P_n) = \sqrt{t} T_\rho^t(P_n^t) = \sqrt{t} Q_n^t = Q_n$

So we get for every $n$: $V_\rho^t(P_n) = (T_\rho^t)^{-1} \circ T_\rho(P_n)$, that validates the composition formula.

$$(3.4) \quad \boxed{V_\rho^t = (T_\rho^t)^{-1} \circ T_\rho}$$

**Inversion of the operator $V_\rho^t$.**

From (3.4) we can write: $(V_\rho^t)^{-1} = (T_\rho)^{-1} \circ T_\rho^t$

Recall now that the definition of equi-normal measure is symmetric, and that the secondary measure associated to the density $\rho_t$ is $t\mu$. So, by changing the roles, $\rho \leftrightarrow \rho_t$, the primary density $\rho$ becomes equi-normal with $\rho_t$ and this with an inverse proportionality coefficient.

So, with the change $t \leftrightarrow \dfrac{1}{t}$ we get: (3.5) $\boxed{(V_\rho^t)^{-1}(f(x)) = \dfrac{1}{t} f(x) + (1-\dfrac{1}{t})(x-c_1)T_\rho^t(f(x))}$

Note : the new parameter $\dfrac{1}{t}$ no longer belongs to the interval $]0,1]$. This is not necessary here because $\rho$ deduced from $\rho_t$ with this coefficient is by primary hypothesis a density of probability.

Application to solving specific integral equations.

We consider here the equation : $\boxed{(E_\lambda): f(x) + \lambda(x-c_1)\int_I \dfrac{f(u)-f(x)}{u-x}\rho(u)du = g(x)}$

with $\lambda > 0$, $g$ a given function and $f$ unknown function in the hyperplane $H_\rho$.

If we note: $t = \dfrac{1}{1+\lambda}$, the equation above easily translates into :

$tf(x) + (1-t)(x-c_1)T_\rho(f)(x) = tg(x)$, or else: $V_\rho^t(f(x)) = tg(x)$, with $t \in ]0,1[$.

According to the precedent results, it can be resolved if $g$ is an element of $H_\rho^t$, using the formula : $f(x) = \dfrac{1}{t}[tg(x)) + (1-\dfrac{1}{t})(x-c_1)T_\rho^t(tg)(x)]$, that leads after simplifications to :

$f(x) = g(x) - \dfrac{\lambda}{1+\lambda}(x-c_1)T_\rho^{\frac{1}{1+\lambda}}(g)(x)$. In integral terms, this solution is made explicit by:

$$(3.6) \quad f(x) = g(x) - \dfrac{\lambda}{1+\lambda}(x-c_1)\int_I \dfrac{g(u) - g(x)}{u-x} \rho_{\frac{1}{1+\lambda}}(u)du$$

Note that the condition $\lambda > 0$ is not necessary, it is sufficient that $\rho_{\frac{1}{1+\lambda}}$ is a probability density function. Note also that the equation ($E_\lambda$) is realized when $g$ is a constant function $g(x) = C$ with the evident solution $f(x) = g(x) = C$. So, thanks linearity, the formula (3.6) can be used with the simple hypothesis: $g$ is element of $L^2(I, \rho_{\frac{1}{1+\lambda}})$.

For example if we chooses the Tchebychev measure of the second kind over $]-1,1[$

$\rho(x) = \dfrac{2}{\pi}\sqrt{1-x^2}$, and $\lambda = \dfrac{-1}{2}$ we have $c_1 = 0, t = 2$ and so $\rho_2(x) = \dfrac{1}{\pi\sqrt{1-x^2}}$ is the Tchebychev density of the first kind.

Then in this case the equation ($E_{-\frac{1}{2}}$): $f(x) - \dfrac{x}{\pi}\int_{-1}^{1} \dfrac{f(t) - f(x)}{t-x}\sqrt{1-t^2}\,dt = g(x)$ has for solution: $f(x) = g(x) + \dfrac{x}{\pi}\int_{-1}^{1} \dfrac{g(t) - g(x)}{t-x} \times \dfrac{dt}{\sqrt{1-t^2}}$

Here bellows some checks using MAPLE:

```
> rho:=x->2*sqrt(1-x^2)/Pi:
> rho2:=x->1/(Pi*sqrt(1-x^2)):
> V:=proc(f)
local t,x,g;
g:=f(x)-(x/2)*int((f(t)-f(x))*rho(t)/(t-x),t=-1..1);
g:=simplify(g);g:=unapply(g,x);end:

> W:=proc(g)
local t,x,f;
f:=g(x)+x*int((g(t)-g(x))*rho2(t)/(t-x),t=-
1..1);f:=simplify(f);f:=unapply(f,x);end:

> V(W(x->2*x^11-7*x^10+8*x^5-3*x+2));
```
$x \to 2x^{11} - 7x^{10} + 8x^5 - 3x + 2$

```
> V(W(x->1/(1+x^2)));
```
$x \to \dfrac{1}{x^2+1}$

```
> V(W(x->x^3/(x+2)));
```
$x \to \dfrac{x^3}{x+2}$

```
> V(W(x->1/(x+3)^2));
```
$x \to \dfrac{1}{x^2+6x+9}$

## 4. Formulas of compositions.

Recall that the secondary measure of $\rho_t$ equals $t\mu$. So, if $t$ and $s$ are two elements defining equi-normal densities $\rho_t$ and $\rho_s$, we have the obvious relation $(\rho_t)_s = \rho_{t.s}$

We also saw in the previous paragraph that the operator $\dfrac{1}{\sqrt{t}} V_\rho^t$ transforms $(P_n)_{n \geq 1}$ orthonormal for the density $\rho$ into $(P_n^t)_{n \geq 1}$ orthonormal for the density $\rho_t$.

From these two points we deduce easily the relation: $\dfrac{1}{\sqrt{t.s}} V_\rho^{t.s} = (\dfrac{1}{\sqrt{s}} V_{\rho_t}^s) \circ (\dfrac{1}{\sqrt{t}} V_\rho^t)$ and so we get the composition formula: (4.1) $\boxed{V_\rho^{t.s} = V_{\rho_t}^s \circ V_\rho^t}$

By changing the initial density into $\rho_u$, the formula above extends to $V_{\rho_u}^{t.s} = V_{\rho_{t.u}}^s \circ V_{\rho_u}^t$

Using the definition: $f(x) \mapsto V_\rho^t(f(x)) = tf(x) + (1-t)(x-c_1)T_\rho(f(x))$, the development of the formula (4.1) leads after simplifications to:

$$(4.2) \quad \boxed{T_\rho^t \circ T_\rho((x-c_1)f(x)) = \left(\dfrac{T_\rho - tT_\rho^t}{1-t}\right)(f(x))}$$

By obvious changing we get the more general formula:

$$(4.3) \quad \boxed{T_\rho^t \circ T_\rho^s((x-c_1)f(x)) = \left(\dfrac{sT_\rho^s - tT_\rho^t}{s-t}\right)(f(x))}$$

If we note for simplifying $S$ the operator: $f(x) \mapsto (x-c_1)f(x)$, we get the barycentric formula above:

$$(4.4) \quad \boxed{T_\rho^t \circ T_\rho^s \circ S = \dfrac{sT_\rho^s - tT_\rho^t}{s-t}}$$

Here is a check using MAPLE in the case of Tchebychev's measures mentioned above:

```
>rho:=x->2*sqrt(1-x^2)/Pi:
>rho2:=x->1/(Pi*sqrt(1-x^2)):

>T:=proc(f)
local t,x,g;
g:=int((f(t)-f(x))*rho(t)/(t-x),t=-1..1);
g:=simplify(g);g:=unapply(g,x);end:

>T2:=proc(f)
local t,x,g;
g:=int((f(t)-f(x))*rho2(t)/(t-x),t=-1..1);
g:=simplify(g);g:=unapply(g,x);end:

>f:=x->7*x^5-4*x^3+x/(x^2+3);
```
$$f := x \to 7x^5 - 4x^3 + \frac{x}{x^2+3}$$

```
>g:=x->x*f(x):
>a:=T2(f)(x):
>b:=T(f)(x):
>simplify(2*a-b);
```
$$\frac{41x^2 - 24\sqrt{3} + 81 + 56x^6 + 178x^4}{8(x^2+3)}$$

```
>T2(T(g))(x);
```
$$-\frac{-41x^2 - 178x^4 - 81 - 56x^6 + 24\sqrt{3}}{8(x^2+3)}$$

Note finally that if we apply the formula (4.3) to a function of type $x \mapsto f(x) = \dfrac{1}{x-z}$ who is an eigenvector for any operator $T_\rho^t$, we obtain the relation linking the Stieltjes transforms $S_\rho^t(z)$ of these equi-normal measures: (4.5) $\boxed{(z-c_1)S_\rho^t(z)S_\rho^s(z) = \dfrac{tS_\rho^t(z) - sS_\rho^s(z)}{t-s}}$

References.

[1] Christian Berg, Annales de la Faculté des Sciences de Toulouse, Sér 6 Vol. S5 (1996), p.9-32.
[2] G.A Baker, Jr and P. Graves-Morris, 'Padé approximants', Cambridge University Press, London 1996.
[3] R. Groux. C.R Acad.Sci. Paris, Ser.I. 2007. Vol 345. pages 373-376.